\documentclass[a4paper,10pt]{article}

\usepackage[breaklinks]{hyperref}
\usepackage{amssymb}
\usepackage{amsthm}
\usepackage{amsmath}
\usepackage{graphicx}
\usepackage{xcolor}

\usepackage{sectsty}
\sectionfont{\rmfamily\mdseries\Large}
\subsectionfont{\rmfamily\mdseries\itshape\large}
\subsubsectionfont{\sffamily\mdseries\slshape}

\newcommand{\figdir}{.}

\usepackage[sc]{mathpazo}

\usepackage{textcomp }

\usepackage[T1]{fontenc} 

\newtheorem{thm}{Theorem}
\newtheorem*{cor}{Corollary}
\newtheorem{lemma}[thm]{Lemma}
\theoremstyle{definition}
\newtheorem*{defn}{Definition}
\newcommand{\R}{\mathbb{R}}

\newcommand{\dd}{\text{\rm d}}
\newcommand{\dvol}{\text{\rm dvol}}

\DeclareMathOperator{\Vol}{vol}
\DeclareMathOperator{\Area}{area}
\DeclareMathOperator{\TSC}{tsc}
\DeclareMathOperator{\Top}{Top}
\allowdisplaybreaks

\title{On the magnitude of spheres, surfaces and other homogeneous spaces}
\author{Simon Willerton}
\date{} 

\begin{document}

\maketitle

\begin{abstract}
In this paper we calculate the magnitude of metric spaces using measures rather than finite 
subsets as had been done previously.  An 
explicit formula for the magnitude of an $n$-sphere with its intrinsic metric is given.  For 
an arbitrary homogeneous Riemannian manifold the leading terms of the asymptotic 
expansion of the magnitude are calculated and expressed in terms of the volume and total 
scalar curvature of the manifold.
\end{abstract}

\tableofcontents

\section*{Introduction}
\addcontentsline{toc}{section}{Introduction}
The magnitude of finite metric spaces is a partially defined numerical invariant which was 
introduced by Leinster~\cite{Leinster:Magnitude} via a 
category theoretic
analogue of the Euler characteristic of spaces, group and posets.   If a finite metric space is scaled up sufficiently then the magnitude tends upwards to the number of points in the metric space and, typically, as the metric space is scaled down the magnitude tends to one; this allows one to think of the magnitude as an `effective number of points' of the metric space. Interestingly, magnitude has applications in measuring biodiversity (see \cite{SolowPolasky:MeasuringBiologicalDiversity} and \cite{Leinster:MaximumEntropy}).

In this paper we look at calculating the magnitude of infinite metric spaces using a measure theoretic approach.  In 
\cite{LeinsterWillerton:AsymptoticMagnitude} and 
\cite{Willerton:Calculations} we
speculatively studied the magnitude of {infinite} metric spaces using 
finite subsets as approximations, then in the first version of this paper we showed that in some specific examples we got the same results by using a measure theoretic approach.  Inspired by this, Mark Meckes~\cite{Meckes:PositiveDefinite} proved that if a metric space is so-called `positive definite' then the two approaches to defining magnitude of the space coincide, provided that the space has a weight measure (see below).  The class of positive definite metric spaces includes spheres with the intrinsic metric and all subsets of Euclidean spaces.

In~\cite{LeinsterWillerton:AsymptoticMagnitude} it is shown using the finite approximation approach that a closed, bounded subinterval of the real line has a magnitude equal to one plus half its length; similarly magnitudes are calculated for circles and Cantor sets.  In Section~1 below, it is shown that these spaces all have weight measures which can straightforwardly be calculated, giving an easy way to obtain the magnitude.

In~\cite{Willerton:Calculations} the magnitudes of various simple subsets 
of $\R^2$ and $\R^3$ were calculated numerically and for various spaces considered it 
was observed that if
the metric space was scaled up then the magnitude seemed to 
asymptotically depend on
classical intrinsic volumes of the space such as volume and Euler 
characteristic.  For \emph{convex} sets it appeared that this dependence held even non-asymptotically, so for a square it seems, numerically at least, that the magnitude is given by $\tfrac{1}{2\pi}\text{area}+\tfrac{1}{4}\text{perimeter}+1$.

%
In this paper the main focus is on the magnitude of homogeneous Riemannian 
manifolds, using the
measure theoretic approach.  The magnitude is calculated explicitly for the $n$-sphere of radius $R$ with its intrinsic metric, leading to the following answer.
\[\left|S^n_R\right|=\begin{cases}
\dfrac{2 \Big(\big(\frac{R}{n-1}\big)^2 + 
1\Big)\Big(\big(\frac{R}{n-3}\big)^2 + 1\Big) \dots
\Bigl(\bigl(\frac{R}{1}\bigr)^2 + 1\Bigr) }{1+e^{-\pi R}}& \text{for 
$n$ even}\\[1.5em]
\dfrac{\pi R \Bigl(\bigl(\frac{R}{n-1}\bigr)^2 + 
1\Bigr)\Bigl(\bigl(\frac{R}{n-3}\bigr)^2 + 1\Bigr)
\dots \Bigl(\bigl(\frac{R}{2}\bigr)^2 + 1\Bigr) }{1-e^{-\pi R}}\quad& 
\text{for $n$ odd}.
\end{cases}
\]
For a general homogeneous Riemannian manifold $X$ we calculate the leading 
three terms in
the magnitude as the manifold is scaled up: the first term is 
proportional to the volume, the
second term vanishes and the third term is proportional to the total 
scalar curvature ($\TSC$).  Writing $tX$ for $X$ scaled up by a factor of $t>0$ and writing $\omega_n$ for the volume of the unit $n$-ball, we have
   \[
       \left|tX\right|
       =
       \tfrac{1}{n!\,\omega_n}\Vol(X)t^n+0\cdot t^{n-1}+
           \tfrac{n+1}{6\,n!\,\omega_n}\TSC(X)t^{n-2}
           +O(t^{n-3})
           \quad\text{as }t\to\infty.
   \]
  This
suffices to determine the asymptotics of the magnitude of a homogeneous 
surface --- such a surface is necessarily a sphere, a projective plane, a torus or a Klein bottle --- the
asymptotic expression involves just the area and Euler characteristic 
of the surface  $\Sigma$:
\[\left|t\Sigma\right|=\tfrac{1}{2\pi}\Area(\Sigma)t^2 + \chi(\Sigma) + 
O(t^{-2})\qquad\text{as }t\to\infty.\]
This expression in the area and Euler characteristic is precisely the  valuation of \cite{Willerton:Calculations}.

There are a few things to note.
The first thing to note is that Riemannian manifolds are not in general positive definite~\cite[Section~3.1]{Meckes:PositiveDefinite} (although spheres \emph{are} positive definite) so the measure theoretic approach is not known to agree with the finite approximation approach in general.  However, the formula for the magnitude of homogeneous Riemannian manifolds agrees with that for the `spread'~\cite{Willerton:spread} and the methods developed here can be applied more generally to the spread of non-homogeneous Riemannian manifolds.

The second thing to note is that in the cases that are known, namely the $n$-spheres and 
surfaces, the constant
         term is precisely the Euler characteristic.  This is 
particularly interesting as
         the magnitude arose from a category theoretic generalization of 
the notion of Euler
         characteristic of various mathematical objects.  Indeed we 
might have named
         magnitude as ``the Euler chacteristic of metric spaces'' if 
that weren't such a
         potentially confusing term.

The final thing to note is that the intrinsic volumes, which are the terms so far identified in the 
asymptotic expansion
         all appear in the asymptotic expansion of the heat kernel
         \cite{Donnelly:HeatEquationTubes} and, as the definition of magnitude is 
slightly
         reminiscent of the heat kernel, this looks intriguing.

Before giving a more detailed synopsis of the contents of the paper, it is worth making 
some comments about the metric that we are 
considering on a
Riemannian manifold.  There is potential for confusion here in that we 
are considering metric spaces and Riemannian metrics, so the word 
`metric' is being used with two different senses.  To reduce the 
likelihood of confusion, in this paragraph and the next paragraph we will talk about 
`distance 
metric' for the metric on a metric space and `Riemannian metric' for an inner product on 
each tangent space of a smooth manifold: a distance metric is a global object, a 
Riemannian metric 
is an infinitessimal object.  In order to define the magnitude of a 
Riemannian manifold we need a distance metric on the manifold, we will 
be considering the intrinsic distance metric which arises in the 
following way.  A Riemannian metric gives rise to an arc-length function 
and this gives rise to a distance metric by defining the distance 
between two points to be the infimum of the lengths of the paths between 
the two points.  

If there is a Riemannian embedding of a Riemannian 
manifold in Euclidean space then there are two obvious distance metrics, the 
intrinsic metric and the subspace metric: in this paper we are 
primarily considering the intrinsic metric whereas in 
\cite{LeinsterWillerton:AsymptoticMagnitude} and 
\cite{Willerton:Calculations} we considered mainly the subspace metric. 
   So, for instance, on the Earth, the intrinsic distance metric is 
realized by travelling over the surface of the Earth, but the subspace 
distance metric is realized by tunnelling \emph{through} the Earth.
We compare the asymptotics for the magnitude of the $n$-sphere with the subspace metric and to that with the intrinsic metric; they are seen to have the same leading term but different lower order terms.   

There now follows a more detailed synopsis of what is in this paper.

\subsection*{Synopsis}
In Section~1 we consider the definition of the magnitude of a metric space using a measure 
on the space.  The key observation here is that if $X$ is a homogeneous 
metric space and $\mu$ is a non-zero homogeneous measure on $X$, then 
the magnitude of $X$ is given by the following formula: for any $y\in X$,
     \[|X|= \frac{\int_{x\in X}\dd\mu}{\int_{x\in X}e^{-d(x,y)}\dd\mu},\]
provided that the denominator is finite and non-zero.  It is 
also shown in this section that the measure theoretic approach can be used to calculate the magnitudes of  the circle, 
line segment and Cantor set, which were calculated using finite approximations  in~\cite{LeinsterWillerton:AsymptoticMagnitude}.

In Section~2 we give a brief exposition on the notion of intrinsic volumes of 
Riemannian manifolds.  Intrinsic volumes are defined for a wide class of 
subsets of Euclidean space and include such invariants as the volume, 
the surface area of the boundary and the Euler characteristic.  In 
\cite{LeinsterWillerton:AsymptoticMagnitude} and 
\cite{Willerton:Calculations} we mainly considered these for polyconvex 
sets --- that is unions of convex sets --- following Klain and Rota's 
treatment \cite{KlainRota:GeometricProbability}, but here we consider a 
differential treatment, utilizing \cite{Morvan:GeneralizedCurvatures} 
and \cite{Gray:Tubes}.  One key point for us is that for $M^n$ a closed 
submanifold of $\R^k$ these intrinsic volumes are intrinsic in the sense 
of differential geometry --- that is to say they are independent of the 
embedding and depend only on the Riemannian metric of $M^n$.  So for 
such a Riemannian manifold there are $n+1$ intrinsic volumes 
$\{\mu_i\}_{i=0}^n$, but it turns out that about half of these vanish 
because $M^n$ is closed, $\mu_i(M^n)=0$ if $n-i$ is odd: the fact 
that $\mu_{n-1}(M^n)=0$ can be explained by the fact that $\mu_{n-1}$ 
measures the surface of the boundary, which is of course trivial when 
$M^n$ is closed.  The non-zero intrinsic volumes are also known in this 
case as Killing-Lipschitz curvatures and in a few cases are well-known 
invariants: $\mu_n(M^n)$ is the volume of $M^n$; the next non-trivial 
intrinsic volume $\mu_{n-2}(M^n)$ is the total scalar curvature (tsc) of 
$M^n$; and $\mu_0(M^n)$ is the Euler characteristic of $M^n$.  In 
general these can be obtained by integrating over $M^n$ a certain form 
defined in terms of the curvature.

In Section~3 we get to one of the main results which is an explicit 
formula for  $\left|S^n_R\right|$ the magnitude of the $n$-sphere of radius $R$ equipped with its intrinsic metric.
(Here it is slightly misleading to talk of the radius of the $n$-sphere 
as it is being thought of as an abstract manifold with the arclength 
distance metric and not as a subset of $\R^{n+1}$ with the subspace 
distance metric; however it is perhaps slightly less comfortable to talk 
of the sphere of curvature $R^{-2}$.)  The formula has many of 
the properties you might expect, such as $|S^n_R|\to 1$ as $R\to 0$. 
Asymptotically as $R\to \infty$, writing $\omega_n$ for the volume 
of the unit $n$-ball, we get
   \begin{align*}
       \left|S^n_R\right|
       &=
       \tfrac{1}{n!\,\omega_n}\Vol(S^n_1)R^n+0\cdot R^{n-1}+
           \tfrac{n+1}{6\,n!\,\omega_n}\TSC(S^n_1)R^{n-2}
           +0\cdot R^{n-3}\\
           &\qquad
             +\dots
           +\chi(S^n_1)
       +O(e^{-R}).
   \end{align*}
This is certainly consistent with the asymptotics being described by the 
intrinsic volumes: for instance the alternate terms vanish and we pick 
up the volume, the total scalar curvature and the Euler characteristic.

In Section~4 we see the other main result that for all homogeneous 
Riemannian manifolds the first three terms in the asymptotic expansion 
are like those in the above expansion, further work will be needed to 
get more terms in the expansion.  However, this is sufficient to show that for a homogeneous 
Riemannian surface  asymptotically the magnitude has the form conjectured for the asymptotic magnitude of 
metric subsets of Euclidean space in 
\cite{LeinsterWillerton:AsymptoticMagnitude} and 
\cite{Willerton:Calculations}.

In Section~5 the magnitude of the $n$-sphere with its subspace metric is considered.  A closed form is given for $n=2$, but not for general $n$;  however, for general $n$ it is shown that the leading term in the asymptotic expansion is the same as for the intrinsic metric, namely $\Vol(S^n_R)/n!\,\omega_n$, but that the next non-trivial term is different.

\subsection*{Acknowlegements}
It is a pleasure to thank Jonathan Jordan, Tom Leinster and David Speyer 
for their contributions to Section~1, as well as for other informative discussions.  I would also like to thank Mark Meckes for useful comments and conversations.

\section{Extending magnitude to infinite metric spaces}
In this section we first recall the definition of the magnitude of finite 
metric spaces in terms of weights and generalize it to infinite metric spaces using weight measures.  We show that Speyer's theorem for homogeneous spaces immediately generalizes  and calculate the magnitude using weight measures for intervals, circles and Cantor sets.

\subsection{Generalizing magnitude using measures}
Leinster~\cite{Leinster:EulerChar} introduced the notion of the Euler characteristic of a finite 
category which is not
always defined, but which generalizes many pre-existing notions of Euler 
characteristic,
such as the Euler characteristic of a group or a poset or the 
cardinality of a set.  Then, using
Lawvere's observation that a metric space can be viewed as an enriched 
category, he
defined~\cite{Leinster:Magnitude} the notion of the Euler characteristic or cardinality of a 
finite metric space --- it was
decided that these terms were already too overburdened, so we renamed 
the concept the
'magnitude' of a metric space.  This has the following definition.
\begin{defn}
For $X$ a finite metric space, a \emph{weighting} on $X$ is an 
assignment of a number
$w_x\in \R$ to each point $x\in X$ such that the \emph{weight equation} 
is satisfied for
each point $y\in X$:
  \[\sum_{x\in X} e^{-d(x,y)}w_x =1.\]
  If a weighting exists for $X$ then the magnitude $|X|$ is defined to 
be the sum of the
weights:
  \[|X|:=\sum_x w_x.\]
  \end{defn}
Note the following two things about the definition.
Firstly, the weights can be negative.
Secondly, if more than one weighting exists then the magnitude is 
independent of the choice of
weighting (a proof appears below in the more general measure theoretic 
setting).

If you require a little intuition then one way of thinking about the 
weight equations is as follows.  Each point can be
thought of as an organism, with the weight being the amount of  `heat' 
being generated by
that point (a negative weight denoting some sort of refrigeration). 
Each organism
experiences the heat generated by the other organisms in a fashion which 
decays
exponentially with the distance.  The organisms each wish to experience 
unit heat and this
is what is expressed by the weight equation.

We can try to extend the idea of magnitude to non-finite metric spaces. 
  Two approaches
present themselves.
\begin{itemize}
\item Use finite approximations to the metric space.
\item Use a measure in the definition of a weighting.
\end{itemize}
The first of these approaches was explored 
in~\cite{LeinsterWillerton:AsymptoticMagnitude} 
and~\cite{Willerton:Calculations}.  In this approach a metric space $X$ is
approximated by taking a sequence $\{X^{(i)}\}_{i=1}^\infty$ of finite 
subsets of $X$, which
converges to $X$ in the Hausdorff topology, so $X^{(i)}\to X$ as 
$i\to\infty$.  One can then
take the corresponding sequence of magnitudes 
$\{|X^{(i)}|\}_{i=1}^\infty$ and hope that
they tend to a limit, and furthermore hope that the limit is independent 
of the choice of approximating sequence.
Following these papers, Mark Meckes~\cite{Meckes:PositiveDefinite} proved that for a large class of spaces, positive definite spaces, which includes subsets of Euclidean spaces, this method gives a definition of magnitude which is independent of the chosen sequence.  This
approach has the advantage that it can be implemented on a computer to 
obtain
approximations to the magnitude of some spaces as in~\cite{Willerton:Calculations}.

The second approach is the one that will be used in this paper.  Inspired by an earlier version of this paper, Mark Meckes showed that for positive definite spaces when this approach can be used, that is `when a weight measure exists', then it is equivalent to the first approach using an approximating sequence, so magnitude can be calculated either way.  This approach  
was independently suggested to me by Tom Leinster and David Speyer (my 
own thoughts on the matter lay more in realm of currents and 
differential forms).  The idea is to generalize
the definition of magnitude using measures in the following way.
\begin{defn}
If $X$ is a metric space then a \emph{weight measure} on $X$ is a finite signed Borel 
measure $\nu$
on $X$ such that for all $y\in X$
\[\int_{x\in X} e^{-d(x,y)}\dd\nu(x)=1.\]
If a weight measure $\nu$ exists then $|X|$, the \emph{magnitude} of $X$, is 
defined to be the mass
of $X$ with respect to $\nu$:
\[|X|:=\int_X\dd\nu.\]
\end{defn}
Note that the magnitude, if it exists, is independent of the choice of 
weight measure.
Suppose that $\nu$ and $\overline{\nu}$ are weight measures then
\begin{align*}
\int_{x\in X}\dd\nu(x)&=\int_{x\in X}\left(\int_{x'\in 
X}e^{-d(x',x)}\dd\overline{\nu}(x')\right)\dd
\nu(x)\\
   &=\int_{x'\in X}\left(\int_{x\in 
X}e^{-d(x,x')}\dd\nu(x)\right)\dd\overline{\nu}(x')
   =\int_{x'\in X}\dd\overline{\nu}(x').
   \end{align*}

It is worth giving a brief mention here to some connections between 
measure theory and
Riemannian geometry (for more details see, for 
example,~\cite{Morvan:GeneralizedCurvatures}).  If $M^n$ is an 
orientable Riemannian manifold then there is a
volume form $\dvol$, well-defined up to sign --- the choice of sign 
corresponding to the
choice of orientation.  The volume form is defined by smoothly picking 
an ordered
orthonormal basis $\partial/\partial x_{1},\dots,\partial/\partial 
x_{n}$ for each tangent
space, taking the dual basis $\dd x_1,\dots, \dd x_n$ and defining 
$\dvol:=\dd x_1\wedge
\dots\wedge \dd x_n$.   Integrating with respect to this volume form 
does not quite give a
measure because of the sign ambiguity --- think how on the real line we 
have $\int_a^b f\,
\dd x=-\int^a_b f\, \dd x$.  However, there is an associated measure 
$\Vol$, called the
Riemannian measure, which obviates the sign problem.  In the case of Euclidean 
space this
measure is the usual Lebesgue measure.  As no confusion should occur, we 
will not
distinguish between the volume form and the associated measure.

If $M^n$ is an orientable $n$-dimensional submanifold of Euclidean space 
$\R^k$ then
there are two obvious measures that one might consider.  Firstly the 
Riemannian metric on
$\R^k$ can be pulled back to $M^n$ to make it into an orientable 
Riemannian manifold so
we can use the associated density.  Alternatively, one could use the 
$n$-dimensional
Hausdorff measure on $\R^k$ restricted to $M^n$.  Provided that the 
Hausdorff measure
is normalized so that the $n$-dimensional cube $[0,1]^n$ has unit 
measure, these two
measures on $M^n$ coincide.

One of the key tools of this paper is the following generalization of 
Speyer's theorem on the magnitude of finite homogeneous spaces given in 
\cite{LeinsterWillerton:AsymptoticMagnitude}.  Recall that a homogeneous 
space is one for which there is a transitive group of isometries, or, in 
other words, one in which all of the points are `the same'.
\begin{thm}[Speyer's Homogeneous Magnitude Theorem]
\label{Thm:Speyer}
Suppose that $X$ is a homogeneous metric space
and $\mu$ is an invariant
measure on $X$ then $\int_{x\in X} e^{-d(x,y)}\dd\mu$ is independent of 
$y$.  If this quantity
is non-zero and finite then a weight measure on
$X$ is given by
\[\frac{\mu}{\int_{x\in X}e^{-d(x,y)}\dd\mu}.\]
Thus the magnitude is given by
\[|X|=\frac{\int_{x\in X}\dd\mu}{\int_{x\in X}e^{-d(x,y)}\dd\mu}.\]
\end{thm}
\begin{proof}
     This follows immediately from the definition of magnitude.
\end{proof}
It is worth pointing out here that every \emph{compact} homogeneous space has a unique-up-to-scaling, invariant, positive measure, thus has a well-defined magnitude.

Aside from the calculations for the interval and Cantor set given below, 
this theorem is the only way I currently know for calculating the 
magnitude of non-finite metric spaces.

\subsection{Agreement with previous calculations}
We can now look at the spaces whose magnitude was defined using finite 
approximations
in~\cite{LeinsterWillerton:AsymptoticMagnitude}, and show that they all have weight 
measures.  These
spaces are closed straight line segments, circles and Cantor sets.

We consider first $L_\ell$ the length $\ell$ closed line segment.  The following was 
essentially told to me by Tom Leinster and David Speyer independently.
\begin{thm}
\label{Thm:IntervalWeightMeasure}
Let $\mu$ be the Lebesgue measure on $L_\ell$ the length $\ell$ line 
segment and let
  $\delta_a$ and $\delta_b$ be the Dirac delta measures supported at  the two respective end
points.  Then a weight measure on $L_\ell$ is given by
   \(\tfrac12 (\delta_a +\delta_b+\mu)\).
  Hence the magnitude is simply:
  \[|L_\ell|=1+\ell/2.\]
\end{thm}
\begin{proof}
It is just necessary to show that the given measure satisfies the weight equation for 
every point $y\in L_\ell$.  We can consider $L_\ell$ as the subset $[a,b]\subset \R$ where $b-a=\ell$, then for every $y\in L_\ell$
\begin{align*}
\int_{x\in 
L_\ell}e^{-d(x,y)}\tfrac12&\bigl(\dd\delta_a+\dd\delta_b+\dd\mu\bigr)\\*
  &= \tfrac12\left(e^{-d(a,y)}+e^{-d(b,y)}+\int_a^y e^{-(y-x)}\dd 
x+\int_y^b e^{-(x-y)}\dd x
\right)\\
  &= \tfrac12\left(e^{a-y}+e^{y-b}+[e^{x-y}]^y_a-[e^{y-x}]_y^b\right)
     = 1.
\end{align*}
Thus the magnitude is given by
\[|L_\ell|=\int_{L_\ell}\tfrac 12 
\bigl(\dd\delta_a+\dd\delta_b+\dd\mu\bigr)=\tfrac12 (1+1+
\ell)=1+\ell/2,\]
as required.
\end{proof}
Now turn to the ternary Cantor set.  The following lemma, pointed out to me by Jonathan Jordan, will be rather 
useful in calculating
the magnitude.  It says that if we have a weight measure on a subset of the line which on some subinterval is a multiple of the Lebesgue measure then if the subinterval is removed then a weight measure is obtained for the resulting space by putting simple Dirac delta measures at the end points of the removed subinterval.
\begin{lemma}
\label{Lemma:MeasureOnHole}
If $X$ is a subset of $\R$ with $\nu$ a weight measure on $X$ with 
$[a,b]\subset X$ and $
\nu|_{[a,b]}=\frac12\mu|_{[a,b]}$, then a weight measure on $X\setminus 
(a,b)$ is given by
\[\nu|_{X\setminus (a,b)}+\tfrac 
12\tanh\bigl(\tfrac{b-a}2\bigr)(\delta_a+\delta_b).\]
So the magnitude satisfies
\[\left|X\setminus (a,b)\right|=|X|-\tfrac{b-a}2 + 
\tanh\bigl(\tfrac{b-a}2\bigr).\]
\end{lemma}
\begin{proof}
We need the easy-to-verify identity
\(\frac{e^A-e^B}{e^A+e^B}=\tanh\bigl(\frac{A-B}2\bigr)\)
to obtain 
$e^{y-a}-e^{y-b}=\tanh\bigl(\frac{b-a}2\bigr)\big(e^{y-a}+e^{y-b}\bigr)$. Thus 
for
$y\le a$
\begin{align*}
\int_{x\in X\setminus 
(a,b)}e^{-d(x,y)}&\bigl[\dd\nu+\tfrac12\tanh\bigl(\tfrac{b-a}2\bigr)(\dd
\delta_a+\dd\delta_b)\bigr]\\
  &=
   \int_{x\in X}e^{-d(x,y)}\dd\nu-\int_{x\in 
[a,b]}e^{-d(x,y)}\tfrac12\dd\mu\\
   &\qquad+\tfrac12\tanh\bigl(\tfrac{b-a}2\bigr)
\int_{x\in X\setminus (a,b)}e^{-d(x,y)}(\dd\delta_a+\dd\delta_b)\\
&=
1+\tfrac12[e^{y-b}-e^{y-a}]+\tfrac12\tanh\bigl(\tfrac{b-a}2\bigr)\bigl(e^{y-a}+e^{y-b}\bigr)\\
&=1.
\end{align*}
The expression for the magnitude follows immediately.
\end{proof}
This can now be repeatedly to obtain the magnitude of the Cantor set, 
giving, as we expect, 
the same answer as in \cite{LeinsterWillerton:AsymptoticMagnitude}.  Note that 
there the Cantor set was
approximated below via finite subsets, and here it is approximated above 
by subsets of
$[0,\ell]$ with non-zero measure.
\begin{thm}
The length $\ell$ ternary Cantor set $T_\ell$ has magnitude given by
\[|T_\ell|=1+\sum_{i=1}^\infty 2^{i-1} 
\tanh\biggl(\frac{\ell}{2.3^i}\biggr).\]
\end{thm}
\begin{proof}
\newcommand{\TT}[1]{T_{\ell,#1}}
\newcommand{\nunu}[1]{\nu_{\ell,#1}}
Start by defining $\TT{0}$ to be a length $\ell$ line interval.
As in Theorem~\ref{Thm:IntervalWeightMeasure}, the interval $\TT{0}$ is given the standard weight measure $\frac{1}{2}(\mu+\delta_0+\delta_\ell)$, which will be called $\nunu{0}$.  Define $\TT{j}$ and its weight measure $\nunu{j}$ inductively by first removing from $\TT{j-1}$ the open middle third of each maximal subinterval; this means removing the point $e\ell$ if $e$ satisfies the following inequality of ternary expansions:
    \[0.a_1a_2\cdots a_{j-1}1 <e <0.a_1a_2\cdots a_{j-1}2 \quad    
        \text{where }a_1,\dots,a_{j-1}\in\{0,2\}.\]
The length of each removed open interval will be $(1/3)^j\ell$.  By Lemma~\ref{Lemma:MeasureOnHole} above we obtain a weight measure $\nunu{j}$ by restricting the weight measure $\nunu{j-1}$ on $\TT{j-1}$, except on the newly pruned end points where we have to put a Dirac delta measure of mass $\frac{1}{2}\tanh((1/3)^j\ell/2)$.

We can say explicitly what the measure is.  Let $S_j$ be the set of newly exposed end points at level $j$, then
  \[\nunu{j}=\frac12\left\{
                          \left.\mu\right|_{\TT{j}}
                           +\delta_0+\delta_{\ell} 
                           +\sum_{i=1}^j
                                    \tanh\biggl(\frac{\ell}{2.3^i}\biggr)
                                    \sum_{s\in S_i} \delta_s 
                         \right\}.
  \]
This gives us a sequence of measures $\{\nunu{j}\}_{j=1}^\infty$ on the interval $[0,\ell]$.  In the Banach space of finite signed Borel measures on $[0,\ell]$, this sequence converges in the total variation norm: $\nunu{j}\to \nu_\ell$ where
  \[\nu_\ell=\frac12\left\{
                      \delta_0+\delta_{\ell} 
                           +\sum_{i=1}^\infty
                                    \tanh\biggl(\frac{\ell}{2.3^i}\biggr)
                                    \sum_{s\in S_i} \delta_s 
                         \right\}.
  \]
Convergence in the total variation norm implies convergence pointwise on any continuous function on $[0,\ell]$, such as $e^{-d({\cdot},y)}$ for $y\in T_\ell$.  Hence
  \[\int_{x\in [0,\ell]} e^{-d(x,y)}\dd\nunu{j} \to \int_{x\in [0,\ell]} e^{-d(x,y)}\dd\nu_\ell.\]
  But as $y\in T_\ell\subset \TT{j}$ for all $j$ and as each $\nunu{j}$ is a weight measure, each of the integrals on the left hand side equals $1$, thus the right hand side is $1$ and $\nu_\ell$ is a weight measure on $T_\ell$.  We can thus calculate the magnitude of $T_\ell$:
 \[|T_\ell| = \int_{T_\ell} \dd \nu_\ell=1+\sum_{i=1}^\infty 2^{i}\frac{1}{2} 
\tanh\biggl(\frac{\ell}{2.3^i}\biggr),\]
as required.
\end{proof}

Now finish this section by considering the circle.
\begin{thm}
The cardinality of $C_\ell$ the circle of circumference $\ell$ is given by
\[|C_\ell|=\frac{\ell}{\int_0^{2\pi}e^{-d(x,y)}\dd x}.\]
\end{thm}
\begin{proof}
This follows immediately from the homogeneity theorem.
\end{proof}
In conclusion, the calculations of magnitude using weight measures are different in these cases to the calculations using finite approximations and could be viewed as being easier.  However, although weight measures exist in these cases, it is not known that weight measures exist for all compact metric spaces, or indeed for all compact subsets of Euclidean space.

\section{Intrinsic volumes for Riemannian manifolds}
\label{Section:IntrinsicVolumes}
This section is a brief exposition on intrinsic volumes for Riemannian manifolds; these 
invariants are also known as curvature measures, Lipschitz-Killing curvatures, generalized 
curvatures and quermassintegrale.  I will mention Weyl's Tube Formula and the connection 
with total scalar curvature and volumes of geodesic speres.  This treatment is mainly 
obtained from the books~\cite{Gray:Tubes} and~\cite{Morvan:GeneralizedCurvatures}.
This section also serves to fix some notation and conventions.

In~\cite{LeinsterWillerton:AsymptoticMagnitude}, following~
\cite{KlainRota:GeometricProbability}, we considered the intrinsic volumes $\{\mu_i
\}_{i=1}^N$ as 
functions on
polyconvex subsets of Euclidean space $\R^N$.  In fact the intrinsic 
volumes are defined
for a much larger class of subsets of $\R^N$ via the theory of curvature 
measures (see
\cite{Morvan:GeneralizedCurvatures}).  The intrinsic volumes can be 
defined for subsets
which have an appropriate notion of unit normal bundle --- a so-called 
normal cycle.  This
class of subsets includes smooth submanifolds with boundary and more 
generally sets of
positive reach.    In general, the
basic properties of the intrinsic volumes are the following.
\begin{description}
\item[Additivity] $\mu_i(X_1\cup X_2)=\mu_i(X_1)+\mu_i(X_2)-\mu_i(X_1\cap X_2)$.
\item[Homogeneity] $\mu_i(tX)=t^i \mu_i(X)$.
\item[Normalization] $\mu_i([0,1]^i)=1$.
\end{description}
Then $\mu_N(X)$ is the usual $N$-dimensional volume of $X$ and 
$\mu_0(X)$ is the Euler
characteristic of $X$.

For a closed smooth submanifold of $\R^N$ there are many formulas for 
the intrinsic
volumes.  The simplest formulas are for  when the submanifold in 
question is  a
hypersurface $H^n\subset\R^{n+1}$. In this case, at each point in $H$ 
there are $n$
principal curvatures $\kappa_1,\dots,\kappa_n$, and the $j$th 
symmetrized curvature
$s_j$ is defined to be the $j$th-elementary symmetric function in the 
principal curvatures: $
\prod_l(1+t\kappa_l)=\sum_js_j t^j$.  The intrinsic volumes are given, 
up to a factor, by
integrating these symmetrized curvatures over the submanifold.  More precisely, 
recalling that $
\sigma_j$ denotes the volume of the unit $j$-sphere, for $0\le i\le n$,
\[\mu_i(H^n)=\begin{cases}
\dfrac{2}{\sigma_{n-i}}\displaystyle\int_H s_{n-i} \,\dvol&\text{if 
$n-i$ is even},
\\
0 &\text{if $n-i$ is odd}.
\end{cases}
\]
 From this one obtains, for instance, $\mu_n(H^n)=\Vol(H^n)$.  This allows us to 
calculate the intrinsic
volumes of $\overline{S}^n_R\subset \R^{n+1}$ the $n$-sphere of radius 
$R$ sitting in $
\R^{n+1}$ in the usual way.  Here all of the principal curvatures are the reciprocal of the 
radius, 
$R^{-1}$, so the
symmetrized curvatures are given by $s_j=\binom{n}{j}R^{-j}$ and so, for 
$0\le i\le n$,
\[\mu_i(\overline{S}_R^n):=\begin{cases}
\dfrac{2\sigma_n}{\sigma_{n-i}}\displaystyle\binom{n}{i}R^i&\text{if 
$n-i$ is even},
\\
0 &\text{if $n-i$ is odd}.
\end{cases}
\]

  In the case of a submanifold of codimension greater than $1$ there are 
not principal
curvatures, but the intrinsic volumes can be expressed in terms of the 
marginally more
complicated Lipschitz-Killing curvatures.  However it is still true that 
for $X^n$ a closed
submanifold $\mu_i(X^n)=0$ whenever $n-i$ is odd.  From these curvature 
expressions, as
Weyl notes, a good undergraduate could derive the following Tube Formula.  As an easy 
exercise, you may wish to verify this for the $n$-sphere using the formulas given above for 
the 
intrinsic volumes.
\begin{thm}[Weyl's Tube Formula]
Suppose that $X^n\subset \R^N$ is a closed, smooth submanifold of 
Euclidean space,
that $B_\epsilon X\subset \R^N$ is the set of points of distance at most 
$\epsilon$ from $X
$, and that $\omega_i$ is the volume of the unit $i$-ball, then
\[\Vol_N(B_\epsilon X)=\sum_{i=0}^N \mu_{N-i}(X)\omega_i \epsilon^i.\]
\end{thm}
Weyl's key observation, however, was that the intrinsic volumes of $X$ 
can be expressed
purely in terms of its \emph{intrinsic} geometry --- we have already 
observed that $
\mu_n(X^n)=\Vol_n(X)$ and $\mu_0(X^n)=\chi(X^n)$.  So the intrinsic 
volumes only depend
on the Riemannian metric of $X$ and not on the way it is embedded, nor 
on the
codimension in which it is embedded.  For instance, the intrinsic 
volumes can be written as
polynomials in the Riemann curvature tensor.  However, in this paper there is 
just one of these
intrinsic volumes that we are interested in.   After the
$n$-volume, the next non-trivial 
intrinsic volume is proportional to the  \emph{total scalar curvature}:
\[\mu_{n-2}(X^n)=\frac{1}{4\pi}\int_{x\in X} \tau(x) \,\dvol.\]
Here $\tau(x)$ is the \emph{scalar curvature} which can be thought of as 
measuring the
``lack of stuff'' near $x$.  For instance, if you draw a circle on the 
surface of the Earth then
it has a smaller circumference than a Euclidean circle with the same radius --- 
this deficiency in the
circumference is measured by the scalar curvature (in the case of a 
surface the scalar
curvature is twice the Gaussian curvature).  To make this precise, for $x\in X$ and $r>0$, 
let 
$F_r(x)\subset X$ be the
set of points at a distance of $r$ from $x$; for $r$ sufficiently small 
this will be an $(n-1)$-sphere, known as a \textit{geodesic sphere}.  The scalar curvature $
\tau(x)$ at $x$ is uniquely determined by the second order term in the power series 
expansion of the volume of geodesic spheres centered at $x$, in the following sense:
  \[\Vol_{n-1}\bigl(F_r(x)\bigr)=\sigma_{n-1}r^{n-1}\left(1-\frac{\tau(x)}{6n}r^2 
+O(r^4)\right)
    \qquad\text{as }r\to 0.\]
It is in this way that the scalar curvature can be viewed as a measure of the ``lack of stuff'' 
near $x$.
  The above formula will be used later in the calculation of sub-leading term in 
the asymptotics of the
magnitude of homogeneous metric spaces.  As an example, note that in an 
$n$-sphere of
radius $R$ we have, for all $x$,
$\Vol_{n-1}\bigl(F_r(x)\bigr)=\sigma_{n-1}\left(R\sin(r/R)\right)^{n-1}$.  So the scalar 
curvature is given by $
\tau=n(n-1)R^{-2}$ which you can check is consistent with the value 
calculated for $
\mu_{n-2}(S^n_R)$ above.

\section{Calculation of the magnitude of an $n$-sphere}
This section consists of one of the main results of this paper which is an explicit formula for the magnitude of the $n$-sphere.  The asymptotic behaviour is also considered and the Euler characteristic is seen to be the constant term in the asymptotic expansion.  

We will use a weight measure to calculate the magnitude, by Theorem~3.6 of~\cite{Meckes:PositiveDefinite} we know that the sphere with its intrinsic metric is positive definite, and so by Theorems~2.3 and~2.4 and Corollary~2.7 of~\cite{Meckes:PositiveDefinite} we know that we would have got the same result using finite approximations.

We take $S_R^n$ the $n$-sphere with the intrinsic metric of curvature 
${1}/{R^2}$
(in other words, having radius $R$) 
thought of as a
homogeneous space.  From Speyer's Homogeneous Magnitude Theorem 
(Theorem~\ref{Thm:Speyer}) the magnitude is given by
\[\left|S^n_R\right|=\frac{\int_{x\in S^n_R} \dvol}{\int_{x\in S^n_R} 
e^{-d(x_0,x)} \dvol} \qquad
\text{for any }x_0\in S^n_R.\]
This can be calculated explicitly.
\begin{thm}
\label{Thm:MagnitudeOfSpheres}
For $n\ge 1$, the magnitude of the homogeneous $n$-sphere of radius $R$ 
  is given by
\[\left|S^n_R\right|=\begin{cases}
\dfrac{2 \Big(\big(\frac{R}{n-1}\big)^2 + 
1\Big)\Big(\big(\frac{R}{n-3}\big)^2 + 1\Big) \dots
\Bigl(\bigl(\frac{R}{1}\bigr)^2 + 1\Bigr) }{1+e^{-\pi R}}& \text{for 
$n$ even}\\[1.5em]
\dfrac{\pi R \Bigl(\bigl(\frac{R}{n-1}\bigr)^2 + 
1\Bigr)\Bigl(\bigl(\frac{R}{n-3}\bigr)^2 + 1\Bigr)
\dots \Bigl(\bigl(\frac{R}{2}\bigr)^2 + 1\Bigr) }{1-e^{-\pi R}}\quad& 
\text{for $n$ odd},
\end{cases}
\]
where for $n=1$ the numerator is interpreted as being $\pi R$.
\end{thm}
\begin{proof}
The idea is to find some explicit integrals to calculate magnitude, 
evaluate them to obtain
$|S^{1}_R |$ and $|S^{2}_R |$, show that the induction formula holds,
   $$|S^{n+2}_R |=\Bigl(\big(\tfrac{R}{n+1}\big)^2+1\Bigr) |S^{n}_R |,$$
and then apply induction.

We first fix some point $x_0$ and define the coordinate $s$ on the 
$n$-sphere to be the
distance from the point $x_0$, so $s(x):=d(x,x_0)$.  The points at a 
distance $s$ from
$x_0$ form an $(n-1)$-sphere $F_{s}$ of Euclidean radius $R\sin(\frac 
sR)$, see Figure~\ref{Fig:LocalBall};
\begin{figure}[tb]
   \[\raisebox{-.45\height}{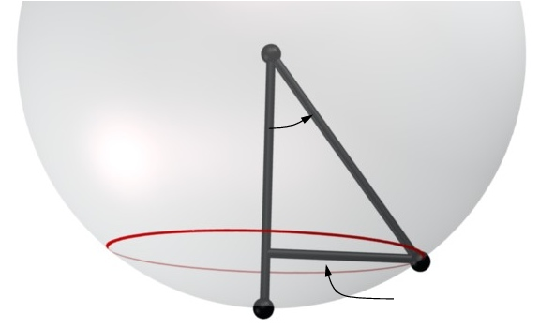}\]
   \caption{The Euclidean radius of the points $F_s$ at a distance
     $s$ from $x_0$ is $R\sin\left(\frac{s}{R}\right)$.}
   \label{Fig:LocalBall}
\end{figure}
thus, writing $
\sigma_{n-1}$ for the volume of the unit $(n-1)$-sphere, the volume of 
$F_{s}$ is given by $
\sigma_{n-1} \bigl(R\sin(\frac{s}{R})\bigr)^{n-1}$. Locally, as $s$ 
parametrizes distance, the volume form on
$S^n_R$ is the product of  $\dd s$ and the volume form on $F_{s}$.  So
\begin{align*}
  \left|S^n_R\right|
    &=
    \frac{\int_{x\in S^n_R} \dvol}{\int_{x\in S^n_R} e^{-d(x_0,x)} \dvol}
    =
    \frac{\int_{s=0}^{\pi R}\int_{F_s} \dd\text{vol}_{F_s}\,\dd s}
           {\int_{s=0}^{\pi R}\int_{F_s} e^{-s}\,\dd\text{vol}_{F_s}\, 
\dd s} \\
   & =
    \frac{\int_{s=0}^{\pi R}\sigma_{n-1} 
\bigl(R\sin(\frac{s}{R})\bigr)^{n-1}\,\dd s}
           {\int_{s=0}^{\pi R}e^{-s}\sigma_{n-1} 
\bigl(R\sin(\frac{s}{R})\bigr)^{n-1}\, \dd s}\\
   &=
    \frac{\int_{r=0}^{\pi}\sin^{n-1}(r)\,\dd r}
           {\int_{r=0}^{\pi }e^{-rR}\sin^{n-1}(r)\, \dd r},
\end{align*}
where we made the substitution $s=rR$.   This can be evaluated for the 
$1$-sphere and
$2$-sphere:
\begin{gather*}
  \left|S^1_R\right|
     =
     \frac{\int_{0}^{\pi}\dd r}
           {\int_{0}^{\pi }e^{-rR}\, \dd r}
     = \frac{\pi}{-\frac{1}{R}[e^{-rR}]_0^\pi}
     = \frac{\pi R}{1-e^{-\pi R}}\;;\\
  \left|S^2_R\right|
     =
     \frac{\int_{0}^{\pi}\sin(r)\,\dd r}
           {\int_{0}^{\pi }e^{-rR}\sin(r)\, \dd r}
     = \frac{2}{\text{Im}\int_{0}^{\pi }e^{-rR}e^{ir}\, \dd r}
    = \frac{2(R^2+1)}{1+e^{-\pi R}}\;.
\end{gather*}
These will form the base cases for the induction.  For the inductive 
step, write
\[K_n:=\int_{0}^{\pi}\sin^{n-1}(r)\,\dd r
\quad\text{and}\quad
I_n:=\int_{0}^{\pi }e^{-rR}\sin^{n-1}(r)\, \dd r
\]
so that  $\left|S^n_R\right|=\frac{K_n}{I_n}$.  We will show that
\[(n+1)K_{n+2}=nK_n\quad\text{and}\quad
(n+1)\Bigl(\big(\tfrac{R}{n+1}\big)^2+1\Bigr)I_{n+2}=nI_n\]
  from which it is immediate that
$|S^{n+2}_R |=\Bigl(\big(\frac{R}{n+1}\big)^2+1\Bigr) |S^{n}_R |$, and 
the required result
follows by induction.

To prove the relation for $K_n$ we observe
\begin{align*}
K_{n+2}
    &=
    \int_{0}^{\pi}\sin^{n+1}(r)\,\dd r
    =
    \int_{0}^{\pi}\left(1-\cos^2(r)\right)\sin^{n-1}(r)\,\dd r\\
    &=
    \int_{0}^{\pi}\sin^{n-1}(r)\,\dd r
    -\int_{0}^{\pi}\cos(r)\left(\cos(r)\sin^{n-1}(r)\right)\dd r\\
     &=
    \int_{0}^{\pi}\sin^{n-1}(r)\,\dd r
-\left[\cos(r)\tfrac1n\sin^{n}(r)\right]^\pi_0+\int_{0}^{\pi}\left(-\sin(r)\right)\left(\tfrac1n
\sin^{n}(r)\right)\dd r\\
    &=K_n-0-\tfrac1n K_{n+2},
\end{align*}
from which we deduce $(n+1)K_{n+2}=nK_n$, as we wanted.

To prove the relation for $I_n$ we perform integration by parts twice:
\begin{align*}
   I_{n+2}
   &=
   \int_{0}^{\pi }e^{-rR}\sin^{n+1}(r)\, \dd r\\
   &=
   \left[-\tfrac 1R e^{-rR}\sin^{n+1}(r)\right]^\pi_0
   -\int_{0}^{\pi }-\tfrac 1Re^{-rR}(n+1)\cos(r)\sin^{n}(r)\, \dd r\\
   &=
   0+
   \tfrac{(n+1)}{R}\left[-\tfrac 1R e^{-rR}\cos(r)\sin^{n}(r)\right]^\pi_0\\
     &\qquad -\tfrac{(n+1)}{R}\int_{0}^{\pi }-\tfrac 
1Re^{-rR}\left(-\sin(r)\sin^{n}(r)
     +n\cos(r)\cos(r)\sin^{n-1}(r)\right) \dd r\\
   &=
   0+0
   +\tfrac{(n+1)}{R^2}\int_{0}^{\pi }e^{-rR}\left(-\sin^{n+1}(r)
        +n\left(1-\sin^2(r)\right)\sin^{n-1}(r)\right) \dd r\\
   &=-\tfrac{(n+1)^2}{R^2}I_{n+2} +\tfrac{(n+1)n}{R^2}I_n.
\end{align*}
 From this we deduce that 
$(n+1)\bigl(\big(\tfrac{R}{n+1}\big)^2+1\bigr)I_{n+2}=nI_n$
which suffices to complete the inductive step and hence we complete the 
proof.
\end{proof}

The first thing to note from this is that the formula for $|S^n_R|$ in 
the above theorem also
works for $S^0_R$ provided that this notation is interpreted correctly. 
   If $S^0_R$ is taken
to mean the two points on the equator of $S_R^1$ the radius $R$ circle, 
then $S_R^0$ is
the metric space consisting of two points a distance $\pi R$ apart, and 
then it is immediate
that the magnitude is given by
\[|S_R^0|=\frac{2}{1+e^{-\pi R}}.\]

The next thing to note is that $|S^n_R|\to 1$ as $R\to 0$, so 
effectively, the $n$-sphere
looks more and more like a point as it gets smaller and smaller.

A third thing to note is that the two formulae can be combined into a 
single formula which
clearly demonstrates the roots of the numerator lying at alternate 
integral points on the
imaginary axis:
\[|S^n_R|=\alpha_n
\frac{\left(R-(n-1)i\right) \left(R-(n-3)i\right)\dots 
\bigl(R-(-n+3)i\bigr)\bigl(R-(-n+1)i\bigr)}
{1+e^{-\pi(R+ni)}}\]
where $i=\sqrt{-1}$ and $\alpha_n$ is a constant which ensures that 
$|S^n_R|\to 1$ as $R\to 0$.  

A fourth thing worthy of note is that the denominators very quickly 
converge to $1$.  This
means that the behaviour for $R$ away from zero, say for $R>4$, is 
dominated by the
numerator.  We can define $P(S^n_R)$, for $n\geqslant 0$, to be the 
numerator:
  \[P(S^n_R):=\begin{cases}
{2 \Big(\big(\frac{R}{n-1}\big)^2 + 1\Big)\Big(\big(\frac{R}{n-3}\big)^2 
+ 1\Big) \dots
\Big(\big(\frac{R}{1}\big)^2 + 1\Big) }& \text{for $n$ even},\\
{\pi R \Big(\big(\frac{R}{n-1}\big)^2 + 
1\Big)\Big(\big(\frac{R}{n-3}\big)^2 + 1\Big) \dots
\Big(\big(\frac{R}{2}\big)^2 + 1\Big) }\quad& \text{for $n$ odd},
\end{cases}
\]
where the empty product is interpreted as $1$ for $n=0$ and $n=1$.  This 
polynomial in $R
$ gives us the asymptotic behaviour of $|S^n_R|$ in the sense that 
$|S^n_R|-P(S^n_R)\to
0$ as $R\to\infty$.  Previously for the closure of an open subset of 
Euclidean space the
asymptotics of the magnitude were conjectured to be controlled by 
certain intrinsic
volumes.  In the case of the spheres the leading term and constant term 
can be seen to
behave in the same way.
\begin{thm}
The leading term of $P(S_R^n)$ is proportional to the volume of the 
$n$-sphere $S^n_R$ and the
constant term is equal to its Euler characteristic: writing $\omega_n$ for the volume of the unit $n$-ball,
   \[P(S^n_R)=\frac{\Vol(S^n_1)R^n}{n!\,\omega_n}+\cdots+\chi(S^n_1).\]
\end{thm}
\begin{proof}
To deal with the Euler characteristic first, note that the Euler 
characteristic of a sphere just
depends on the parity of the dimension: $\chi(S^{2i})=2$ and 
$\chi(S^{2i+1})=0$.  The
constant term in the polynomial $P(S^n_R)$ is immediately seen to have 
exactly this
property.

That the leading term has the required form can be proved inductively using the inductive formula 
 \[P(S^n_R)=\bigl(\tfrac{R}{n-1}\bigr)^2P\bigl(S^{n-2}_R\bigr)\]
together with 
 the following relations (see, for example,~\cite{KlainRota:GeometricProbability}) for $\omega_n$ the volume of the unit $n$-ball and $\sigma_n$ the volume of the unit $n$-sphere:
  \[\omega_n=\tfrac{2\pi}{n}\omega_{n-2}
    \quad\text{and}\quad
    \sigma_n=\tfrac{2\pi}{n-1}\sigma_{n-2}.\]
Assuming that the result is true for $P(S^{n-2}_R)$ and writing $\Top$ for the top coefficient of a polynomial, we have
  \begin{align*}
    \Top\Bigl(P\bigl(S^n_R\bigr)\Bigr)
    &=\Top\Bigl(\bigl(\tfrac{R}{n-1}\bigr)^2
        P\bigl(S^{n-2}_R\bigr)\Bigr)
    =\Top\Bigl(\bigl(\tfrac{R}{n-1}\bigr)^2\Bigr)
      \Top\left(P\bigl(S^{n-2}_R\bigr)\right)\\
    &=\frac{1}{(n-1)^2}
       \frac{\Vol\bigl(S^{n-2}_1\bigr)}{(n-2)!\,\omega_{n-2}}
    = \frac{\sigma_{n-2}}{(n-1)^2\,(n-2)!\,\omega_{n-2}}\\
    &= \frac{\frac{n-1}{2\pi}\sigma_{n}}
             {(n-1)(n-1)!\,\frac{n}{2\pi}\omega_{n}}
    =\frac{\sigma_{n}}
             {n!\,\omega_{n}}
    = \frac{\Vol\bigl(S^{n}_1\bigr)}{n!\,\omega_{n}}.
  \end{align*}
The result follows by induction after checking the base cases, $n=1$ and $n=2$, using $\sigma_1=2\pi$, $\omega_1=2$, $\sigma_2=4\pi$ and $\omega_2=\pi$.   
\end{proof}
We will see in the next section that the leading term in the asymptotic expansion of the magnitude of any homogeneous is of this form.  In~\cite{LeinsterWillerton:AsymptoticMagnitude} and~\cite{Willerton:Calculations} we gave a precise conjecture for the asymptotic behaviour of the magnitude of `nice' subsets of Euclidean space and this involved the being asymptotically the same as the `penguin valuation' defined for $X$ a subset of Euclidean space as
  \[P(X):=\sum_i\frac{\mu_i(X)}{i!\,\omega_i}\]
We are here considering the sphere with its intrinsic metric so it is not a metric subspace of Euclidean space, but the leading term and constant term we obtain are the same as for in the penguin valuation.  However, the subdominant term is different.  We would expect ${\mu_{n-2}\bigl(S^n_R\bigr)}/({(n-2)!\,
\omega_{n-2}})$ but we get the  following: it can be proved inductively as the formula for the leading term was in the above theorem, it is also a consequence of the theorem in the next section on general homogeneous Riemannian manifolds.
\begin{thm}
The first non-trivial, subdominant term of $P(S_R^n)$ is given by
\[ \frac{(n+1)}{3(n-1)}\frac{\mu_{n-2}\bigl(S^n_R\bigr)}{(n-2)!\,
\omega_{n-2}}. \]
\end{thm}
These results are compared with analogous ones for the magnitude of the $n$-sphere with the subspace metric in Section~5.

\section{Asymptotics for homogeneous Riemannian manifolds}
In this section we will look at the three leading terms in the 
asymptotics of the magnitude of
a homogeneous Riemannian manifold as the manifold is scaled up.  The 
leading term will
be seen to be proportional to the volume of the manifold, the 
second-leading term will be seen to be zero and the third-leading term 
will be seen to be
proportional to its total scalar curvature, or in other words the top 
two leading terms are
proportional to the top two intrinsic volumes.  In particular, for a 
homogeneous surface $
\Sigma$, the magnitude is asymptotically 
$\frac{1}{2\pi}\Area(\Sigma)+\chi(\Sigma)$, and
this is precisely the evaluation of the functional from 
\cite{LeinsterWillerton:AsymptoticMagnitude}.

We take the magnitude of a Riemannian manifold to mean the weight measure definition; the fact that Riemannian manifolds are not all positive definite~\cite[Section~3.1]{Meckes:PositiveDefinite} means that we do not know this is the same as the finite approximation definition.   However, this still gives an interesting definition, not least because this definition of magnitude for homogeneous Riemannian manifolds agrees with the definition of `spread' from~\cite{Willerton:spread} --- spread is an invariant of metric spaces inspired by biodiversity --- and the methods here generalise to the spread of non-homogeneous Riemannian manifolds.  For details see~\cite{Willerton:spread}.

By Speyer's Homogeneous Magnitude Theorem (Theorem~\ref{Thm:Speyer}) we 
know that the magnitude of a homogeneous
Riemannian manifold $X$ is given by
   \[|X|=\frac{\int_X \dvol}{\int_X e^{-d(x,x_0)}\dvol}.\]
If we look at this as we scale up the space we should consider
\[|t X|=\frac{\int_{t X} \dvol}{\int_{t X} e^{-d(tx,tx_0)}\dvol}
=\frac{t^n\int_ X \dvol}{t^n\int_ X e^{-td(x,x_0)}\dvol}
=\frac{\int_ X \dvol}{\int_ X e^{-td(x,x_0)}\dvol}.\]
We wish to consider the asymptotics of the denominator ${\int_ X 
e^{-td(x,x_0)}\dvol}$.
The relevant piece of asymptotic analysis we will require is a version 
of Watson's Lemma (see 
for example~\cite{Miller:AppliedAsymptoticAnalysis}).
\begin{thm}[Watson's Lemma]
If $g\colon[0,c]\to [0,\infty)$ is a bounded function which, for some 
$N>0$, has an
expansion
\[g(r)=\sum_{i=0}^N\alpha_i r^i+O(r^{N+1}) \qquad\text{as }r\to 0\]
then
\[\int_0^c e^{-tr}g(r)\dd r= \sum_{i=0}^N 
\frac{i!\,\alpha_i}{t^{i+1}}+O(t^{-N-2})\qquad
\text{as }t\to\infty.\]\qed
\end{thm}
The idea behind Watson's Lemma is that as a function of $r\in[0,c]$, 
when $t$ gets
larger and larger $e^{-tr}$ looks more and more like a function 
supported on a smaller and
smaller neighbourhood of $r=0$, so the integral will, up to 
exponentially small corrections,
only depend on the germ of $g$ near $r=0$.

This can now be used to identify the leading asymptotics of the 
magnitude for
homogeneous manifolds.
\begin{thm}
\label{Thm:HomogeneousAsymptotics}
If $X$ is an $n$-dimensional homogeneous Riemannian manifold and $\TSC$ 
denotes total
scalar curvature then as $X$ is scaled up the asymptotics of the 
magnitude are as follows:
  \begin{align*}
   |tX|
    &=
    \frac{1}{n!\,\omega_n}\Bigl(\Vol(tX) + \frac{n+1}{6}\TSC(tX) 
+O(t^{n-4})\Bigr)\\*
    &=
    \frac{\mu_n(X)}{n!\,\omega_n}t^n + \frac{(n+1)\,\mu_{n-2}(X)}{3(n-1)!\,
\omega_{n-2}}t^{n-2} +O(t^{n-4}),\quad \text{as }t\to\infty.
\end{align*}
\end{thm}
\begin{proof}
The idea is to pick a point $x_0\in X$ and define the coordinate 
$r(x)=d(x,x_0)$.  Let $F_r$
be the set of points at a distance $r$ from $x_0$, then if $\tau$ 
denotes the scalar
curvature of $X$, which will be the same at every point, then we know
  \[\Vol_{n-1}(F_r)=\sigma_{n-1}r^{n-1}\Bigl(1-\frac{\tau}{6n}r^2 
+O(r^4)\Bigr)
    \qquad\text{as }r\to 0.\]
Because $r$ parametrizes distance, locally the volume form splits as 
$\dvol_{F_r}\,\dd r$.
We can write $D$ for the diameter of $X$.  Then we find, using Watson's 
Lemma
appropriately,
  \begin{align*}
  \int_ X e^{-td(x,x_0)}\dvol
  &=
  \int_{r=0}^D \int_{F_r}e^{-tr}\dvol_{F_r}\,\dd r
  =
  \int_{r=0}^D e^{-tr}\Vol(F_r)\,\dd r\\
  &=
  \sigma_{n-1}\Bigl(  \frac{(n-1)!}{t^n} -\frac{(n+1)!\,\tau}{6n 
t^{n+2}} + O(t^{-n-4}) \Bigr)
  \quad\text{as }t\to \infty\\
  &= \frac{n!\,\omega_n}{t^n}\Bigl( 1- 
\frac{(n+1)\tau}{6t^2}+O(t^{-4})\Bigr),
  \end{align*}
where in the last equality we have used the fact that 
$\sigma_{n-1}=n\omega_n$.  Thus we
get
\begin{align*}
|tX|&=\frac{\int_ X \dvol}{\int_ X e^{-td(x,x_0)}\dvol}
   =
   \frac{{t^n}\int_ X \dvol}{{n!\,\omega_n}\Bigl( 1- 
\frac{(n+1)\tau}{6t^2}+O(t^{-4})\Bigr)}\\
   &=   \frac{1}{n!\,\omega_n}\Bigl( t^n\int_X\,\dvol + 
\frac{n+1}{6}t^{n-2}\int_X \tau\,\dvol
+O(t^{n-4})\Bigr)\quad\text{as }t\to \infty\\
   &=  \frac{1}{n!\,\omega_n}\Bigl(\Vol(tX) + \frac{n+1}{6}\TSC(tX) 
+O(t^{n-4})\Bigr) .
\end{align*}
This gives the first description of the asymptotics, to get the second 
it is just necessary to
recall that $\mu_{n-2}(X)=\frac{1}{4\pi}\TSC(X)$ that 
$\mu_{n-2}(tX)=t^{n-2}\mu_{n-2}(X)$ and that $\omega_n=\frac{2\pi}{n}
\omega_{n-2}$.
\end{proof}

This theorem immediately gives the asymptotics of the magnitude of a 
homogeneous surface, by taking $n=2$ in the above and noting that 
$\mu_0$ is precisely the Euler characteristic $\chi$.  Note that a homogeneous surface is necessarily either a sphere, a torus, a Klein bottle or a projective plane.
\begin{cor}
For $\Sigma$ a homogeneous Riemannian surface, the magnitude is 
asymptotically given in terms of the area and the Euler characteristic by
\[|t\Sigma|=\frac{\Area(t\Sigma)}{2\pi}+\chi(t\Sigma)+O(t^{-2})\qquad 
\text{as }t\to\infty.\]
\end{cor}
It is interesting to note that in the examples we know of, namely for 
spheres and surfaces, the constant term in the asymptotic expansion of 
the magnitude is precisely the Euler characteristic.

%
%
%

\section{The subspace metric on the sphere}
In this section we look at the metric on the $n$-sphere obtained by thinking of it as sub-metric space of $\R^{n+1}$ via the standard embedding.  An explicit formula is given for the magnitude of the $2$-sphere with this metric, whereas for $n>2$ just the leading few terms of the asymptotic expansion are given.  This is sufficient to see that the leading term is the same as for the magnitude of the sphere with the intrinsic metric but that the lower order terms are not the same.

In \cite{LeinsterWillerton:AsymptoticMagnitude} we showed that equipping the radius $R$ \emph{circle} with the intrinsic metric and the subspace metric gave rise to different magnitudes, but that asymptotically the magnitudes differed by a vanishingly small piece:
  \[\bigl|S_R^1\bigr|,\bigl|S_{\text{sub},R}^1\bigr| = \pi R +O(R^{-1}) 
   \quad\text{as }R\to\infty.\]
  With the techniques so far developed in this paper it is not difficult to calculate the magnitude of the $2$-sphere with these two metrics and see that only the leading term in the asymptotic expansion are the same.
\begin{thm}
For the $2$-sphere of radius $R$ with the subspace metric the magnitude is given by
  \begin{align*} \bigl|S_{{\rm sub},R}^2\bigr| &= \frac{2 R^2}{1-e^{-2R}(1+2R)}\\
      &=2R^2 +O(R^{-1})\quad\text{as }R\to\infty.
\intertext{
For the $2$-sphere of radius $R$ with the intrinsic metric the magnitude is given by} \bigl|S_{R}^2\bigr| &= \frac{2 R^2+2}{1-e^{-\pi R}}\\
      &=2R^2 +2+O(R^{-1})\quad\text{as }R\to\infty.
  \end{align*}
Thus the sub-leading terms of the asymptotic behaviour are different for the two different metrics.
\end{thm}
\begin{proof}The intrinsic metric result is in Theorem~\ref{Thm:MagnitudeOfSpheres}, so it is just necessary to do the calculation for the subspace metric.  This is analogous to the calculation for the intrinsic metric.  By the Homogeneous Magnitude Theorem
\begin{align*} 
    \bigl|S_{{\rm sub},R}^2\bigr|&=    
        \frac{\int_{x\in S^2_R} \dvol}{\int_{x\in S^2_R} e^{-d_{{\rm sub}}(x_0,x)} \dvol}
        =
        \frac{4\pi R^2}{\int_{\theta=0}^\pi e^{-R2\sin(\theta/2)}2\pi R\sin(\theta) R\,\dd\theta}\\
        &=
        \frac{2}{\int_{\theta=0}^\pi e^{-R2\sin(\theta/2)}\sin(\theta)\, \dd\theta}.        
\end{align*}
Using the substitution $s:=2\sin(\theta/2)$ gives the following evaluation of the integral:
  \[\int_{\theta=0}^\pi e^{-R2\sin(\theta/2)}\sin(\theta)\, \dd\theta
    =
    \int_{s=0}^2 e^{-Rs}s\, \dd s
    =
    \frac{1-e^{-2R}(1+2R)}{R^2}.\]
The result follows immediately.
\end{proof}

Whilst I don't know how to get a general closed form for the magnitude of the sphere with the subspace metric, using the techniques of the previous section it is straight forward to calculate the first few terms in the asymptotic expansion of the magnitude and compare them with those for the intrinsic metric.
\begin{thm}
For the $n$-sphere of radius $R$ with the subspace metric the magnitude is given asymptotically by
  \[ \bigl|S_{{\rm sub},R}^n\bigr|=\frac{\Vol(S_R^{n})}{n!\,\omega_n}\biggr(1+\frac{(n+1)n(n-2)}{8}R^{-2}+O(R^{-4})\biggl)\quad\text{as }R\to \infty.\]
For the $n$-sphere of radius $R$ with the intrinsic metric the magnitude is given asymptotically by
  \[ \bigl|S_{R}^n\bigr|=\frac{\Vol(S_R^{n})}{n!\,\omega_n}\biggr(1+\frac{(n+1)n(n-1)}{6}R^{-2}+O(R^{-4})\biggl)\quad\text{as }R\to \infty.\]
 \end{thm}
\begin{proof}The intrinsic metric result follows from Theorem~\ref{Thm:HomogeneousAsymptotics} together with the fact from Section~\ref{Section:IntrinsicVolumes} that the scalar curvature of the $n$-sphere with radius $R$ is $n(n-1)R^{-2}$.

For the subspace metric result, proceed as usual by using Speyer's Homogeneous Magnitude Theorem and writing $\sigma_n$ for the volume of the unit $n$-sphere.
\begin{align*} 
    \bigl|S_{{\rm sub},R}^n\bigr|&=    
        \frac{\int_{x\in S^n_R} \dvol}{\int_{x\in S^n_R} e^{-d_{{\rm sub}}(x_0,x)} \dvol}
        =
        \frac{\Vol(S_R^n)}{\int_{\theta=0}^\pi e^{-R2\sin(\theta/2)}\sigma_{n-1}\bigl(R\sin(\theta)\bigr)^{n-1} R\,\dd\theta}\\
        &=
        \frac{\Vol(S_R^n)}{\sigma_{n-1}R^n\int_{\theta=0}^\pi e^{-R2\sin(\theta/2)}\sin^{n-1}(\theta)\,\dd\theta}.
\end{align*}
Again making the substitution $s:=2\sin(\theta/2)$ whilst noting that $\sin(\theta)\,\dd\theta=s\,\dd s$ and $\sin(\theta)=s\sqrt{1-s^2/4}$ the integral can be manipulated as follows:
  \begin{align*}\int_{\theta=0}^\pi e^{-R2\sin(\theta/2)}&\sin^{n-1}(\theta)\,\dd\theta
  =
  \int_{s=0}^2 e^{-Rs}s^{n-1}(1-{s^2}/{4})^{(n-2)/2}\dd s,
  \intertext{so by the generalized Binomial Theorem,}
  &=
  \int_{s=0}^2 e^{-Rs}\sum_{i=0}^\infty (-1)^i\begin{pmatrix}\tfrac{n-2}{2}\\i\end{pmatrix}
  \frac{s^{n-1+2i}}{4^i}\,\dd s
  \intertext{whence by Watson's Lemma, as $R\to \infty$,}
  &=
  \sum_{i=0}^m (-1)^i\begin{pmatrix}\tfrac{n-2}{2}\\i\end{pmatrix}
  \frac{(n-1+2i)!}{4^iR^{n+2i}}  +O(R^{-n-2m+1})\\
  &=
  \frac{(n-1)!}{R^n}\biggl(1-\frac{(n-2)}{2}\frac{(n+1)!}{(n-1)!\,4}R^{-2}+O(R^{-3})\biggr).
  \end{align*}
 Substituting this back into the above formula gives the required result on recalling that $\sigma_n=n\omega_n$.
\end{proof}

One can think of this result as saying that just as the intrinsic and subspace metrics are infinitesimally the same in a certain sense --- though I don't know the correct way to formalize that --- so the corresponding magnitudes are asymptotically the same in a certain sense, namely $\left|S^n_R\right|/\bigl|S^n_{\text{sub},R}\bigr| \to 1$ as $R\to \infty$, although they are not the same in a stronger sense because, in general, 
$\left|S^n_R\right|-\bigl|S^n_{\text{sub},R}\bigr| \not\to 0$ as $R\to \infty$.


\begin{thebibliography}{99}
\addcontentsline{toc}{section}{References}
\bibitem{Donnelly:HeatEquationTubes}
     H.~Donnelly,
     \textit{Heat equation and the volume of tubes},
     Inventiones Mathematicae (1975) \textbf{29} (3) pp.~239-243.
     %
\bibitem{Gray:Tubes}
     A.~Gray,
     Tubes,
     Birkh\"auser  (2004).
%
\bibitem{KlainRota:GeometricProbability}
     D.~Klain and G.-C.~Rota,
     Introduction to Geometric Probability,
     Cambridge University Press  (1997).
%
\bibitem{Leinster:EulerChar}
    T.~Leinster,
    \textit{The Euler characteristic of a category},
       Documenta Mathematica \textbf{13}  (2008), 21--49.  \url{http://www.math.uni-bielefeld.de/documenta/vol-13/02.html}    
%
\bibitem{Leinster:MaximumEntropy}
    T.~Leinster,
    \textit{A maximum entropy theorem with applications to the measurement of
            biodiversity},
        arXiv preprint. \href{http://arxiv.org/abs/0910.0906}{arxiv:0910.0906}.
%
\bibitem{Leinster:Magnitude} T.~Leinster,
  \textit{The magnitude of metric spaces},
  arXiv preprint. \href{http://arxiv.org/abs/1012.5857v3}{arxiv:1012.5857v3}
%
\bibitem{LeinsterWillerton:AsymptoticMagnitude}
     T.~Leinster and S.~Willerton,
     \textit{On the asymptotic magnitude of metric spaces},
     Geometriae Dedicata, August 2012, 1--24. \href{http://arxiv.org/abs/0908.1582}{arxiv:0908.1582}.
%
\bibitem{Meckes:PositiveDefinite} M.~Meckes,
  \textit{Positive definite metric spaces},
  Positivity, September 2012, 1--25.
 \href{http://arxiv.org/abs/1012.5863}{arxiv:1012.5863}
 %
\bibitem{Miller:AppliedAsymptoticAnalysis} 
    P.~D.~Miller,
   {Applied Asymptotic Analysis},
    Graduate Studies in Mathematics 75, American Mathematical Society,
    2006.
%
\bibitem{Morvan:GeneralizedCurvatures}
     J.-M.~Morvan,
    Generalized Curvatures,
    Springer-Verlag  (2008).
%
\bibitem{Schanuel:Potato}
     S.~Schanuel,
     \textit{What is the length of a potato? An introduction to 
geometric measure theory},
     in {Categories in continuum physics (Buffalo, N.Y., 1982)},   118--126,
     Lecture Notes in Mathematics, \textbf{1174}, Springer, Berlin, 1986.
%
\bibitem{SolowPolasky:MeasuringBiologicalDiversity}
    A.~Solow and S.~Polasky,
    \textit{Measuring biological diversity},
    Environmental and Ecological Statistics (1994) \textbf{1} (2) 95--103.
%
\bibitem{Willerton:Calculations}
     S.~Willerton,
     \textit{Heuristic and computer calculations of the magnitude of 
metric spaces},
   arXiv  preprint. \href{http://arxiv.org/abs/0910.5500}{arxiv:0910.5500}.
%
\bibitem{Willerton:spread} 
     S.~Willerton,
    \textit{Spread: a measure of the size of metric spaces},
     arXiv  preprint. \href{http://arxiv.org/abs/1209.2300}{arXiv:1209.2300}.
%
\bibitem{Weyl:Tubes}
     H.~Weyl,
     \textit{On the volume of tubes},
     American Journal of Mathematics (1939) \textbf{61} (2) pp.~461--472.


\end{thebibliography}
\end{document}